\documentclass[12pt,leqno]{article}
\usepackage{amsmath,amsthm,amsfonts,amssymb,amscd}
\usepackage[latin1]{inputenc}
\usepackage{graphics}
\usepackage{graphicx}

\title{\bf Regular Foliations and Trace Divisors}
\date{\today}

\author{Paulo Sad}

\begin {document}

\maketitle

\begin{abstract}
\noindent A example is given of a  divisor of a curve which is not a trace divisor of a foliation.\footnote{ MSC Number: 32S65}

\end {abstract}

\section {\bf Introduction} 

We deal in this paper with a very concrete question about existence of  foliations with prescribed singularities.

Let $ F=0 $ be the equation of a plane smooth curve $Q$ of degree $4$; we fix the line at infinity $L_{\infty} \subset {\mathbb C}P^2$, transverse to $Q$, and put $D_{\infty} = Q \cap L_{\infty}$. 
We select $16$ points $Q_1, \dots, Q_{16}$ outside $L_{\infty}$ and blow up ${\mathbb C}P^2$ at these points; the curve $Q$ becomes a curve $\hat Q \subset \widehat{ {\mathbb C}P^2}$.  

Given the quartic $Q$ and a foliation $\mathcal F$ which leaves $Q$ invariant we say that the configuration of points $\{Q_1,\dots,Q_{16}\}$ is {\it pre-regular}  for  $\mathcal F$ when $\hat Q$ is a leaf of the foliation $\hat {\mathcal F}$ obtained after blowing up at each of these points (in particular, $\hat Q$ has no singularities of $\hat {\mathcal F}$). A simple example is given by the pencil generated by $Q$ and another transversal quartic; if $\{Q_1,\dots,Q_{16}\}$ is the intersection of the two quadrics, then the divisor  $[ \sum_{j=1}^{16} Q_j -4D_{\infty}]$ is a principal divisor (for simplicity we write $D=0$ for a principal divisor $D$). 

 In (\cite {S}) it is observed that if the points $Q_1,\dots,Q_{16}$ are singularities of $\mathcal F$ then the divisor
$[\sum_{j=1}^{16} m_j Q_j - L D_{\infty}]$ (its $\it{trace\, divisor\, along\, Q}$) defined putting $m_j$ as the tangent multiplicity of $\mathcal F$ along $Q$ at $Q_j$ and $4L= \sum_{j=1}^{16} m_j$ is a principal divisor.  A generic choice of the configuration of points  satisfies none of these "resonance" conditions no matter what the choice of the $m_1,\dots,m_{16}$ is; in particular such a  configuration will not be pre-regular.  We may  then allow resonances for the configuration and ask whether it  is pre-regular for some foliation. This seems to be a difficult problem.  Here in this paper  we give an example of a configuration  satisfying $[\sum_{j=1}^{16} Q_j -4D_{\infty}]\ne 0$ and  $2 [\sum_{j=1}^{16} Q_j -4D_{\infty}]=0$ which  is  pre-regular for no foliation having $2[\sum_{j=1}^{16} Q_j -4D_{\infty}]$ as its trace divisor.

This is related to the so called $\it {Ueda\,\, type}$ of $\hat Q$ (\cite {U}), as studied  (in great generality) in (\cite {P}). When $D=\sum_{j=1}^{16} Q_j -4D_{\infty}$ satisfies $mD=0$ for some integer smaller or equal  to the Ueda type of $\hat Q$, then the configuration is pre-regular for some foliaton. It is interesting to look for  cases where $mD=0$ but the configuration is not pre-regular for any foliation (the simplest case being $2D=0$ and Ueda type of $\hat Q$ equal to 1).

I would like to thank J.V.Pereira and M. Falla Luza for very helpful conversations.

\section {The example}
 
The set $\{Q_1,\dots, Q_{16}\}$ is to be chosen for now as the zero set in $Q$ of the function $G=l_1l_2C^2$, where $C$ is a curve of degree 3 transverse to $Q$ and  $l_1$ and $l_2$ are bitangent lines to $Q$; we assume that the points of $l_1=0$ and $l_2=0$ lying in $Q$ are not contained in $C$. It follows that $D\ne0$ and $2D=0$, for $D=\sum_{j=1}^{16}Q_j$, This example was introduced by Neeman in (\cite {N}). We intend to prove the

\vspace {0.1in}
\noindent {\bf Theorem} Let $Q$ be a generic smooth quartic. There exists a Neeman example such that the configuration of points is not pre-regular for any foliation which has $2D$ as its trace divisor along $Q$.

\vspace {0.1in}
We will fix  $C= F_X$  and assume: 1) the curve $C$ is transverse to $Q$; 2) $F_Y \neq 0$ at the points of $l_1l_2=0$ (there is a third condition that will appear in the proof of the Theorem). The Neeman example mentioned in the statement is the set $\{l_1l_2C^2=0\}\cap Q$.

\vspace {0.1in}
We intend to study a foliation $\mathcal F$ which has $Q$ as an invariant set and which is singular along $Q$  exactly at the points $Q_1,\dots, Q_{16}$; it is assumed that the singularities have multiplicities (along $Q$) equal to 2.  Then  $\mathcal F$ has degree 10  and it is defined by the polynomial 1-form
\begin {equation} \label {def 5}
\Omega ={\tilde G}dF + F{\tilde {\eta}}
\end {equation}

\noindent  with $\tilde {\eta}= BdX -AdY$;  the polynomials $A$ and $B$ have both degree 7.

Since the quotient $\dfrac {{\tilde G}}{G}$ is holomorphic along $Q$, we have that ${\tilde G} = G +\lambda F$ for some $\lambda \in \mathbb C$; therefore 
$$
\Omega = (G+ \lambda F) dF + F{\tilde \eta}= GdF + F(\tilde  \eta + \lambda FdF) = GdF + F{\eta}
$$

\noindent where $\eta = \tilde \eta + \lambda FdF$. This means that we may use $G$ in the expression of $\Omega$.

Let us remark that more generally we may state

\vspace{0.1in}
\noindent  {\bf Proposition 1} Assume $H$ and $\tilde H$ are polynomials such that the curves $H=0$ and $\tilde H=0$ intersect $Q$ (in ${\mathbb C}P^2$) exactly at the points $Q_1,\dots,Q_{16}$ with the same multiplicities.
If $\mathcal H$ is a foliation defined by $HdF+F\xi=0$ which has all its singularities along $Q$  at the points $Q_1,\dots,Q_{16}$, then $\cal H$ is also defined by $\tilde HdF+F{\tilde \xi}$ for some polynomial 1-form $\tilde \xi$.
\vspace{0.1in}

From now on we consider a foliation $\mathcal F$ defined by a polynomial 1-form $\Omega=GdF+F(BdX-AdY)$ which has singularities along $Q$ exactly at the points $Q_1, \dots ,Q_{16}$. We assume that $\mathcal F$ has multiplicity greater or equal to $2$  along $Q$ at any of these points. Let us
now  take local coordinates $(x,y)$ in a neighborhood of some $Q_j$
$$
x=X,  \,\,\,y=F(X,Y)
$$
\noindent or, equivalently
$$
X=x,\,\,\,Y=l(x,y)
$$

\noindent The relations 
$$
F_X(x,l(x,y)) + l_x(x,y).F_Y(x,l(x,y))=0,\,\,\,\, F_Y(x,l(x,y)).l_y(x,y)=1
$$

\noindent follow from $F(x,l(x,y))\equiv y$.

We remark that for each small $c\in \mathbb C$ the map $x\mapsto (x,l(x,c))$ parametrizes the curve $F(X,Y)=c$.

In these local coordinates the 1-form $\Omega$ becomes
\begin {equation} \label {def 6}
g(x,y)dy+y[B(x,l)dx - A(x,l)(l_xdx+l_ydy)]=
\end {equation}
$$
[g-yl_yA(x,l)]dy +y[B(x,l)-l_xA(x,l)]dx, 
$$

\noindent where $g=g(x,y)=G(x,l(x,y))$ and $l=l(x,y),l_x=l_x(x,y), l_y=l_y(x,y)$.

\vspace{0.1 in}
$$
{\dot x}=-g+y\alpha(x,y)),\,\,\,{\dot y}=y\beta(x,y)
$$

\noindent The associated vector field is
$$
[-g(x,y)+yl_yA(x,l(x,y)) \dfrac {\partial}{\partial x} + y[B(x,l(x,y))-l_xA(x,l(x,y))] \dfrac {\partial}{\partial y}
$$
\vspace{0.1in}

\noindent or simply 
$$
{\dot x}=-g(x,y)+y\alpha(x,y) ,\,\,\,{\dot y}=\beta (x,y)
$$

\noindent for  $\alpha(x,y)=l_yA(x,l(x,y))$ and $\beta(x,y)=B(x,l(x,y)) -l_xA(x,l(x,y))$. 

Let us write $g(x,y)=\lambda x^m +y{\tilde g}(x,y)$ and ${\tilde \alpha}(x,y)= \alpha (x,y)-{\tilde g}(x,y)$;  we have then 

\begin {equation} \label {def 7} 
{\dot x}= -\lambda x^m+y{\tilde \alpha}(x,y), \,\,\, {\dot y}=y\beta(x,y)
\end {equation}

\section {Blowing-up the singularities}

In this section we analyse necessary conditions a singularity as in  $(7)$ must satisfy in order to be a $\it pre$-$\it regular$ singularity of the foliation ${\mathcal F}$. 

\vspace {0.1in}
\noindent {\bf Definition 1} A singularity  in $Q$ is ${\it pre}$-${\it regular}$ when after one blow-up the intersection of the exceptional divisor with the strict transform of $Q$ is a regular point of the strict transform of $\mathcal F$.
\vspace {0.1in}

After one blow-up of (7) for the case $m=2$ (putting $y=tx$)  we get
$$
{\dot x}=-\lambda x^2+tx{\tilde \alpha}(x,tx),\,\,\,\,{\dot t}=t\beta(x,tx)-t[-\lambda x + t{\tilde \alpha}(x,tx)];
$$

\noindent let us write ${\tilde \alpha}(x,y)={\tilde \alpha}+ {{\tilde \alpha}_x}x +{{\tilde \alpha}_y}y+ \dots$ and ${ \beta}(x,y)={ \beta}+ {{ \beta}_x}x +{{ \beta}_y}y+ \dots.$

\vspace {0.1in}
 Since the Camacho-Sad index of  a pre-regular singularity is equal to $1$, we have $\beta_x=-\lambda$. The necessary conditions for pre-regularity are then

\begin {equation} \label {def 8}
{\tilde \alpha}=0,\,\,\beta=0,\,\,{\beta}_y-{\tilde \alpha}_x=0,\,\,{\tilde \alpha}_y=0.
\end {equation}

Let us remark that if the singularity was
$$
{\dot x}= -\lambda x^m+y{\tilde \alpha}(x,y), \,\,\, {\dot y}=y\beta(x,y)
$$

\noindent for  $m>2$, then the conditions in (8) of pre-regularity  

$$
{\tilde \alpha}=0,\,\,\beta=0,\,\,{\beta}_y-{\tilde \alpha}_x=0,\,\,{\tilde \alpha}_y=0,\,\, \beta_x=0
$$

\noindent are also, but not all, necessary conditions; the fact that the Camacho-Sad index is $1$ corresponds to ${\beta}^{m-1}_x=-\lambda$, where ${\beta}^{m-1}_x x^{m-1}$ is the $(m-1)$-term in the Taylor development of $\beta(x,y)$.

We proceed to write the conditions in terms of $G,A$ and $B$; for simplicity we assume the singularity to be the point $(0,0)\in \mathbb C^2$.

\vspace {0.1in}
\noindent ${\bf 1)} \,\,  {\tilde \alpha}=0$.  
\vspace{0.1in}

\noindent Since $\tilde \alpha (x,y) = \alpha (x,y) -\tilde g (x,y)$ we have $\alpha (0,0)= \tilde g (0,0)$. But $g(x,y)=\lambda x^2 +y\tilde g (x,y)$; therefore $\dfrac{\partial g}{\partial y}(x,y)= \tilde g (x,y)+y\dfrac {\partial \tilde g}{\partial y}(x,y)$ and $\dfrac{\partial g}{\partial y}(0,0)= \tilde g (0,0)$

\noindent Now we have $\alpha(0,0)=\dfrac{\partial g}{\partial y}(0,0)$ so  $l_y(0,0)A(0,l(0,0))=\dfrac {\partial G}{\partial Y}(0,l(0,0)).l_y(0,0)$.

We conclude that 

\begin {equation} \label {def 9}
 A(0,0)=\dfrac {\partial G}{\partial Y}(0,0)
\end {equation}

\vspace {0.1in}
\noindent ${\bf 2)} \,\,\beta=0$
\vspace {0.1in}

\noindent We have  $B(0,0)-l_x(0,0)A(0,0)=0$ and $l_x(0,0)=-\dfrac{F_X}{F_Y}(0,0)$
\noindent so that 

\begin {equation} \label {def 10}
 A(0,0)F_X(0,0)+ B(0,0)F_Y(0,0)=0.
\end {equation}

\vspace {0.1in}
\noindent ${\bf 3)}\,\, {\tilde \alpha_y}=0$
\vspace {0.1in}

\noindent As before, $\tilde \alpha(x,y)=\alpha (x,y) - \tilde g (x,y)$ and  $\tilde \alpha_y (x,y)= \alpha_y (x,y) - \tilde g_y (x,y)$.

\noindent Since $\alpha (x,y)= l_y(x,y) A(x,l(x,y))$ we get 

$$
\alpha_y(x,y)= l_{yy}(x,y)A(x,l(x,y)) +l_y^2 (x,y)A_Y(x,l(x,y))
$$

 \noindent As before,  $\dfrac{\partial g}{\partial y}(x,y)= \tilde g (x,y)+y\dfrac {\partial \tilde g}{\partial y}(x,y)$; this implies

$$
g_{yy}(x,y)=2\tilde g_y(x,y) +y\tilde g_{yy}(x,y)
$$

\noindent and

$$
l_{yy}(0,0)G_Y(0,0) + l_y^2 (0,0)A_Y(0,0) - \frac {g_{yy}(0,0)}{2} =0
$$

We have now to compute $g_{yy}(0,0)$, $l_y(0,0)$ and $l_{yy}(0,0)$.

\begin {itemize}
\item   $ g(x,y)=G(x,l(x,y))$ and   $g_y(x,y)= G_Y(x,l(x,y))l_y(x,y)$ so 

$g_{yy}(x,y)= l_{yy}(x,y)G_Y(x,l(x,y)) + l_y^2(x,y)G_{YY}(x,l(x,y))$ 

and $l_y^2(0,0)A_Y(0,0)=\dfrac {l_y^2(0,0)G_{YY}(0,0)-  l_{yy}(0,0)G_Y(0,0)}{2}$.

\item  $F(x,l(x,y))=y$ implies  $F_Y(x,l(x,y))l_y(x,y)=1$ 

 and $F_{YY}(x,l(x,y))l_y^2(x,y) + F_Y(x,l(x,y)l_{yy}(x,y)=0$ so that 

$l_{yy}(x,y)=-\dfrac {F_{YY}(x,l(x,y))l_y^2(x,y)}{F_Y}$.
\end {itemize}

Finally we get

\begin {equation} \label {def 12}
A_Y(0,0)= \dfrac{G_{YY}(0,0)}{2} + \dfrac {F_{YY}(0,0)G_Y(0,0)}{2F_Y(0,0)}
\end {equation}

\vspace {0.1in}
\noindent {\bf 4)} $\beta_y-\tilde \alpha_x=0$
\vspace {0.1 in}

\noindent From $\beta (x,y)= B(x,l(x,y))-l_x(x,y)A(x,l(x,y))$ we get 
$$
\beta_y(x,y)=B_Y(x,l(x,y))l_y(x,y)-l_{xy}(x,y)A(x,l(x,y))-l_x(x,y)l_y(x,y)A_Y(x,l(x,y))
$$

\noindent In order to compute $\tilde \alpha_x =\alpha_x - \tilde g_x$ we use:

\begin {itemize}
\item $\alpha(x,y)=l_y(x,y)A(x,l(x,y))$ implies 
$$
\alpha_x (x,y)=l_{xy}(x,y)A(x,l(x,y)) + l_y(x,y[(A_X(x,l(x,y))+ l_x(x,y)A_Y(x,l(x,y))]
$$
\item $g_y(x,y)=\tilde g(x,y)+y\tilde g(x,y)$ implies
$$
\tilde g_x(x,y)= g_{xy}(x,y)-y\tilde g_x(x,y)=
$$
$$
(G_Y . l_y)_x(x,y)-y\tilde g_x(x,y)=
$$
$$
 [G_{YX}(x,l(x,y))+G_{YY}(x,l(x,y))l_x(x,y)]l_y(x,y)+
$$
$$
G_Y(x,l(x,y)).l_{xy}(x,y)-y\tilde g_x(x,y)
$$

\end {itemize}

\noindent It follows that
$$
\tilde \alpha_x(0,0)=l_y(0,0)[A_X(0,0)+l_x(0,0)A_Y(0,0)-G_{YX}(0,0)-G_{YY}l_x(0,0)]
$$

\noindent Using (9) and (11) we have $[B_Y(0,0)-A_X(0,0)]l_y(0,0)=$
$$
l_{xy}(0,0)G_Y(0,0)+ \dfrac {F_{YY}(0,0)G_Y(0,0)}{F_Y(0,0)}l_{xy}(0,0) -G_{XY}(0,0)l_y(0,0)
$$

\noindent To finish the computation we notice that from $l_y(x,y)F_Y(x,l(x,y))=1$ it follows that 
$$
l_{xy}(x,y)F_Y(x,l(x,y))+l_y(x,y)[F_{YX}(x,l(x,y))+F_{YY}(x,l(x,l(x,y))l_x(x,y)=0
$$
\noindent and finally 

\begin {equation}\label {def 12}
 B_Y(0,0)-A_X(0,0)=-\dfrac{F_{XY}(0,0)G_Y(0,0)}{F_Y(0,0)}-G_{XY}(0,0)
\end {equation}

\vspace {0.1in}
\noindent ${\bf 5)}\,\, \beta_x=-\lambda$
\vspace {0.1in}

\noindent We have $\lambda= -\dfrac {g_{xx}(0,0)}{2}$.

\noindent From  $\beta(x,y)=B(x,l(x,y))-l_x(x,y)A(x,l(x,y))$ it follows
$$
\beta_x=B_X(0,0) + [B_Y(0,0)-A_X(0,0)]l_x(0,0)-l_{xx}(0,0)A(0,0)-l_x^2(0,0)  A_Y(0,0)
$$

\noindent On the other hand $g(x,y)=G(x,l(x,y))$ implies 
$$
g_{xx}(0,0)=G_{XX}(0,0) +2G_{XY}(0,0)+l_{xx}(0,0)G_Y(0,0)+G_{YY}(0,0)l_x^2(0,0)
$$

\noindent and $F(x,l(x,y))+ l_x(x,y)F_Y(x,l(x,y))=0$ implies 
$$
l_{xx}(0,0)F_Y(0,0)=-2F_{XY}(0,0)l_x(0,0)- F_{XX}(0,0)-l_x^2(0,0)F_{YY}(0,0)
$$

\noindent Using (11) and (12) we finally get

\begin {equation} \label {def 13}
B_X(0,0)=-\dfrac {G_{XX}(0,0)}{2}- \dfrac {F_{XX}(0,0)G_Y(0,0)}{2F_Y(0,0)}
\end {equation}

\section {Consequences and Proof of the Theorem}

Let us consider again the  foliation $\mathcal F_0$  of degree 10 defined by the 1-form 
 $\Omega =GdF + F( BdX -AdY)$; here $G=l_1l_2C^2$ and $A$, $B$ are polynomials of degree 7. 
The question we address now is: is it possible that all the singularities  $Q_1,\dots,Q_{16}$ are pre-regular? 

We denote by $Q_{13},Q_{14}$ the points of $\{l_1(X,Y)=0\}\cap Q$ and by $Q_{15},Q_{16}$ the points of $\{l_2(X,Y)=0\} \cap Q$; we have $\dfrac{(l_1l_2)_X}{(l_1l_2)_Y}= \dfrac {F_X}{F_Y}$ at these points.

We have the following list of identities

\begin {itemize}
\item $G=(l_1l_2)C^2$
\item $G_X=(l_1l_2)_X C^2 + 2(l_1l_2)CC_X$
\item $G_Y=(l_1l_2)_YC^2 + 2(l_1l_2)CC_Y$
\item $G_{XX}=(l_1l_2)_{XX}C^2 + 4(l_1l_2)_XCC_X+2(l_1l_2)CC_{XX}+2(l_1l_2)C_X^2$
\item $G_{XY}=(l_1l_2)_{XY} C^2+2[(l_1l_2)_X CC_Y+(l_1l_2)_Y C C_X+l_1l_2C_YC_X +l_1l_2 C C_{XY}]$
\item $G_{YY}=(l_1l_2)_{YY}C^2+ 4(l_1l_2)_YCC_Y+2(l_1l_2)C_Y^2+2(l_1l_2)CC)_{YY}$
\end {itemize}

Let us remind that $C=F_X$ and  $F_{XX}\neq 0$ at the points of $\{F_X=0\} \cap Q$.

\vspace {0.1in}
\noindent{\bf Lemma 1} (i) $A=B=0$ at the points of $Q\cap \{C=0\}$;  (ii) $B_X= -(l_1l_2)C_X^2$ at
 $Q_1,\dots,Q_{12}$; (iii) $B_X= -\dfrac {1}{2}(l_1l_2)_{XX}C^2 -2(l_1l_2)_XCC_X - \frac{1}{2}(l_1l_2)_XF_XF_{XX}$ at the point $Q_{13}.$
\proof It is enough to use (9) and (10) and the identities above.
\qed

\vspace{0.1in}

\vspace {0.1in}
\noindent {\bf Proposition 2} \,There exist polynomials $h,h^{\prime}$ of degree 4 and $k,k^{\prime}$  polynomials of degree 3 such that
\begin {equation}
A=hC+kF,\,\,\, B=h^{\prime}C+k^{\prime}F
\end {equation}
\proof   Since the meromorphic function $(\dfrac{A}{C})|_Q$ has its polar divisor supported in $L_{\infty}$, there exists a polynomial $h$ such that $(\dfrac{A}{C})|_Q= h|_Q$. As $(A-hC)|_Q=0$,  it follows that there exists a polynomial
 $k$ satisfying $A-hC=kF$. The same argument works for $B$.   \qed

\vspace {0.1in}
We observe that $(h^{\prime},k^{\prime})$ can be changed by $(h^{\prime}+\mu F,k^{\prime}-\mu C)$ for any $\mu \in {\mathbb C}$. Consequently we may assume that at $Q_{13}$ we have $h^{\prime} \neq0$ and $k^{\prime} \neq 0$.

\vspace {0.1in}
\noindent{\bf Lemma 2}\,\, $h^{\prime}+(l_1l_2)F_{XX}+(l_1l_2)_XF_X=0$ at the points $Q_1,\dots,Q_{16}$.
\proof Let us look first at the points of $\{F_X\} \cap Q=0$. Since $B_X=h^{\prime}F_{XX}$ and $B_X=-(l_1l_2)C_X^2$, we have $h^{\prime}=-(l_1l_2)F_{XX}$. As for $Q_{13},\dots.Q_{16}$, $A=(l_1l_2)_YF_X^2$ (because of (5)). But $B=h^{\prime}F_X$ and $AF_X +BF_Y=0$ (because of (6)), so we get $h^{\prime}=-(l_1l_2)_XF_X$ (using  $\dfrac{(l_1l_2)_X}{(l_1l_2)_Y}= \dfrac {F_X}{F_Y}$).
\qed

\vspace {0.1in}
Let us put $u=h^{\prime}+(l_1l_2)F_{XX}+(l_1l_2)_XF_X$; we wish to show that $u$ is not identically zero along $Q$. In order to do that, we parametrize a neighborhood  of $Q_{13}=(0,0)$ with the variable $X$ and compute the derivative of $u|_Q$. We observe that $(u|_Q)^{\prime}= u_X - \frac{F_X}{F_Y}u_Y$ at this point.

Without any loss of generality we may assume that $l_1(X,Y)=Y+aX$ and $l_2(X,Y)=Y+aX+b$, for some $a\neq 0$ and $b \neq 0$. Using Lemma 1 we get after a straight computation that at $Q_{13}$:
$$
(u|Q)^{\prime}= -a^2F_X  -\frac{1}{2}abF_{XX} -a^2bF_{XY} - (k^{\prime} + ah^{\prime}_Y)
$$

In order to get rid of the term $k^{\prime} + ah^{\prime}_Y$, we replace the 1-form $GdF-F \eta$ that defines the foliation by $G d{\hat F}-{\hat F}\eta$ where $\hat F=cF$ for some $c \in {\mathbb C}$. Since, according to Proposition 2, $B={\hat h}^{\prime}C + {\hat k}^{\prime}{\hat F}$ we see that ${\hat h}^{\prime}=h^{\prime}$ and ${\hat k}^{\prime}=c^{-1}k^{\prime}$. Consequently ${\hat k}^{\prime} +a{\hat h}^{\prime}_Y=c^{-1}k^{\prime}+ah^{\prime}_Y$; we can therefore choose $c\in {\mathbb C}$
in order that ${\hat k}^{\prime} +a{\hat h}^{\prime}_Y=0$. The corresponding  function $\hat u$ that replaces $u$ satisfies

$$
({\hat u}|Q)^{\prime}=c [-a^2F_X  -\frac{1}{2}abF_{XX} -a^2bF_{XY}]
$$

We may now complete the meaning of the term "generic" which appears in the statement of the Theorem.  We ask that the following inequality holds  in $Q_{13}$:

$$
aF_X  +\frac{1}{2}bF_{XX} +abF_{XY}\neq 0 
$$

\noindent This ensures that ${\hat u}|Q$ in not identically zero along $Q$.

\vspace {0.2in}
In order to finish the proof of the Theorem we use the degree 4 polynomial $\hat u$. Since this polynomial  vanishes at $Q_1,\dots,Q_{16}$, we get that the divisor $\sum_{j=1}^{16}Q_j -4D_{\infty}$ is a principal divisor, contradiction.

\vspace {0.3in}
Paulo Sad

{\rm IMPA- Instituto de Matematica Pura e Aplicada}

sad@impa.br

\end{document}